\documentclass[12pt,reqno,psamsfonts]{amsart}
\usepackage{amssymb}

\DeclareMathAlphabet{\Bi}{OT1}{cmm}{b}{it}

\newcommand{\NN}{\mathbb{N}}

\newcommand{\RR}{\mathbb{R}}
\newcommand{\ZZ}{\mathbb{Z}}

\newcommand{\cm}{com\-mu\-ta\-tive mo\-no\-id}
\newcommand{\com}{com\-mu\-ta\-tive or\-der\-ed
mo\-no\-id}
\newcommand{\poag}{par\-tial\-ly or\-der\-ed a\-bel\-ian
group}

\newcommand{\lsc}{lo\-wer se\-mi\-con\-ti\-nu\-ous}
\newcommand{\usc}{up\-per se\-mi\-con\-ti\-nu\-ous}
\newcommand{\anog}{Ar\-chi\-me\-dean norm-com\-plete
di\-men\-sion group with or\-der-u\-nit}
\newcommand{\anogs}{Ar\-chi\-me\-dean norm-com\-plete
di\-men\-sion groups with or\-der-u\-nit}
\newcommand{\grot}{\mathrm{Grot}}

\newcommand{\La}{\boldsymbol{\varLambda}}

\newcommand{\ess}{\boldsymbol{S}}

\newcommand{\bphi}{\boldsymbol{\phi}}

\newcommand{\bpsi}{\boldsymbol{\psi}}

\newcommand{\dnw}{\mathop{\downarrow}\nolimits}
\newcommand{\upw}{\mathop{\uparrow}\nolimits}
\newcommand{\res}{\mathop{\restriction}\nolimits}

\newcommand{\ddnw}
{\mathop{\dnw\kern-8pt\lower2pt\hbox{$\dnw$}}}

\def\a{\mathbf{a}}\def\ag{\mathfrak{a}}
\def\b{\mathbf{b}}\def\bg{\mathfrak{b}}
\def\c{\mathbf{c}}\def\cg{\mathfrak{c}}
\def\d{\mathbf{d}}\def\dg{\mathfrak{d}}

\def\x{\mathbf{x}}
\def\y{\mathbf{y}}

\numberwithin{equation}{section}

\theoremstyle{plain}
\newtheorem{theorem}{Theorem}[section]
\newtheorem{proposition}[theorem]{Proposition}
\newtheorem{lemma}[theorem]{Lemma}
\newtheorem{corollary}[theorem]{Corollary}
\newtheorem{examplepf}[theorem]{Example}
\newtheorem*{ED}{Edwards' Separation Theorem}

\theoremstyle{definition}
\newtheorem{parag}[theorem]{}
\newtheorem{definition}[theorem]{Definition}

\newtheorem{problem}{Problem}
\newtheorem*{ackn}{Acknowledgement}

\theoremstyle{remark}

\begin{document}

\title[Norm-closed intervals]%
{Norm-closed intervals of
norm-complete ordered abelian groups}

\author[F.~Wehrung]{Friedrich~Wehrung}
\address{D\'epartement de Math\'ematiques\\
         Universit\'e de Caen\\
         14032 Caen Cedex\\
         France}
 \email{gremlin@math.unicaen.fr}

 \date{\today}

\begin{abstract} Let $(G,u)$ be an \anog. Continuing a
previous paper, we study
\emph{intervals} (\emph{i.e.}, non\-emp\-ty upward directed
lower subsets) of $G$ which are
\emph{closed} with respect to the canonical norm of
$(G,u)$. In particular, we establish a canonical one-to-one
correspondence between closed intervals of
$G$ and certain affine lower semicontinuous functions on the
state space of $(G,u)$, which allows us to solve several
problems of K. R. Goodearl about inserting affine continuous
functions between convex \usc\ and concave \lsc\ functions.
This yields in turn new results about analogues of
multiplier groups for norm-closed intervals.
\end{abstract}

\maketitle

\section*{Introduction}

A fundamental result about affine continuous functions on
Choquet simplexes is the following one, due to D. A.
Edwards (see \cite[Th\'eor\`eme]{Edw65},
\cite[Theorem II.3.10]{Alf71} or
\cite[Theorem 11.13]{Good86}):

\begin{ED}
Let $K$ be a Choquet simplex and let
$p\colon K\to\RR\cup\{-\infty\}$ and
$q\colon K\to\RR\cup\{+\infty\}$ be functions such that
$p$ is convex \usc, $q$ is concave \lsc\ and $p\leq q$
(componentwise). Then there exists an affine continuous
function
$f\colon K\to\RR$ such that $p\leq f\leq q$.
\end{ED}

If one strengthens the conclusion by requiring the affine
continuous function $f$ to lie in a given subgroup $G$
containing $1$ of the partially ordered abelian group
$\mathrm{Aff}(K)$ of all affine continuous functions on
$K$, then more stringent assumptions on $p$ and $q$ are
necessary --- at least, for all $x\in K$, there should
exist $f\in G$ such that
$p(x)\leq f(x)\leq g(x)$. In some cases, minor variations
around the latter assumption turn out to be sufficient, as
in \cite[Theorem 13.5]{Good86} and \cite[Theorem 3.5]{Pard},
or \cite[Theorem 16.18]{Good86} in the case where $G$
satisfies countable interpolation, with
\cite[Example 15.13]{Good86} to show that it is
not the case that every such ``reasonable" statement
actually holds.\smallskip

Now let $(G,u)$ be an \anog, let $S$ be the state space of
$(G,u)$, let $\phi\colon G\to\mathrm{Aff}(S)$ be the
natural map and let $p\colon S\to\{-\infty\}\cup\RR$ and
$q\colon S\to\RR\cup\{+\infty\}$ be functions such that $p$
is convex \usc, $q$ is concave \lsc\ and $p\leq q$. One
asks whether there exists an element
$x$ of $G$ such that $p\leq\phi(x)\leq q$, under various
additional assumptions on $p$ and $q$. In
\cite[Problem 13]{Good86}, the additional assumption is
that for every discrete extremal state $s$, both $p(s)$ and
$q(s)$ belong to $s[G]\cup\{-\infty,+\infty\}$. In
\cite[Problem 19]{Good86}, the additional assumption is
that $G$ has countable interpolation and for every discrete
extremal state $s$,
$[p(s),\,q(s)]\cap s[G]$ is non empty.\smallskip

We solve both problems here (Theorem~\ref{T:SandwSig} for
Problem 13 and Example~\ref{E:NonSandw} for Problem 19), by
continuing the study, initiated in
\cite{Wehr96}, of monoids of \emph{intervals} (that is, non
empty upward directed lower subsets) of
\poag s. In fact, we will focus here on intervals which are
\emph{closed} with respect to the canonical norm (see
\cite{Good86}) on a \poag\ with order-unit.\smallskip

Furthermore, this study will allow us,
in the ``good" cases, to give an exact
characterization of norm-closed intervals in terms of
affine \lsc\ functions. More specifically, if $(G,u)$ is an
\anog, if $S$ is the state space of $(G,u)$ and if
$\phi\colon G\to\mathrm{Aff}(S)$ is the natural
homomorphism, then to every interval $\ag$ of $G$, one
associates the supremum $q$ of all $\phi(x)$ where
$x\in\ag$. Then $q$ is an affine \lsc\ function from $S$ to
$\RR\cup\{+\infty\}$, and, for every discrete extremal state
$s$ on $S$, $q(s)$ belongs to $s[G]\cup\{+\infty\}$. The
main result of this paper (Theorem~\ref{T:CharSig}) is a
converse of this statement, generalizing to
\emph{norm-closed intervals} the result already known for
\emph{elements} of an \anog, see
\cite[Theorem 15.7]{Good86}.
As an application of this result, analogues of multiplier
groups for \anogs\ with respect to a bounded positive
norm-closed interval are norm-complete
(Theorem~\ref{T:MClArch}).

\section*{Notation and Terminology}

As in \cite{Wehr96}, we will widely use in this paper the
results and notations of \cite{Good86}. Section~\ref{Prel}
will be devoted to prepare the framework of the whole paper.
It recalls in particular some of the ``refinement axioms"
(more specifically, IA, WIA, RD, REF and REF$'$) already
introduced in \cite[Section 1]{Wehr96}.\smallskip

We will denote by $\sqcup$ the disjoint union of sets. If
$X$ is a subset of a set $S$ (understood from the context),
we will denote by $\chi_X$ the characteristic function of
$X$. If $f$ is a function of domain $X$, we will sometimes
use the notation $f=\langle f(x)\colon x\in X\rangle$;
moreover, if $Y$ is a set, we will denote by $f[Y]$ (resp.
$f^{-1}Y$) the direct (resp. inverse) image of $Y$ under
$f$. Following
\cite{Good86}, we will denote by $\ZZ^+$ the set of all
non-negative integers, and put
$\NN=\ZZ^+\setminus\{0\}$.\smallskip

As in \cite{Wehr96}, if
$X$ is a topological space and $M$ is an additive submonoid
of $\RR$, we will denote by $\mathbf{C}(X,M)$ (resp.
$\mathbf{LSC}(X,M)$,
$\mathbf{LSC}_\mathrm{b}(X,M)$) the space of all
real-valued continuous (resp. \lsc, bounded \lsc)
functions from $X$ to
$\RR$; furthermore, if $M$ is an additive submonoid of
$\RR^+$, let $\mathbf{LSC}^\pm(X,M)$ (resp.
$\mathbf{LSC}_\mathrm{b}^\pm(X,M)$) be the ordered additive
subgroup of all differences $f-g$ where both $f$ and $g$
belong to $\mathbf{LSC}(X,M)$ (resp.
$\mathbf{LSC}_\mathrm{b}(X,M)$), with the positive cone
$\mathbf{LSC}(X,M)$ (resp. $\mathbf{LSC}_\mathrm{b}(X,M)$).
\smallskip

We will denote by $\beta\ZZ^+$ the topological space of
all ultrafilters of $\ZZ^+$ (\v Cech-Stone
compactification of the discrete space $\ZZ^+$).\smallskip

If $(P,\leq)$ is a partially ordered set and both $X$ and
$Y$ are subsets of $P$, then we will abbreviate the
statement \((\forall x\in X)(\forall y\in Y)(x\leq y)\) by
\(X\leq Y\). Furthermore, if \(X=\{a_1,\ldots,a_m\}\) and
\(Y=\{b_1,\ldots,b_n\}\), then we will write
\(a_1,\ldots,a_m\le b_1,\ldots,b_n\). If $\alpha$ and
$\beta$ are two cardinal numbers, then we will say that
\emph{$P$ has the
$(\alpha,\beta)$-interpolation property} when for all
\emph{non\-emp\-ty} subsets $X$ and $Y$ of $P$ such that
$|X|\leq\alpha$ and $|Y|\leq\beta$ and $X\leq Y$, there
exists $z\in P$ such that
$X\le\{z\}\leq Y$. The \emph{interpolation property} is the
$(2,2)$-interpolation property. Say that
$P$ is \emph{directed} when for all $x$, $y\in P$, there
exists $z\in P$ such that $x,y\leq z$.

If $X$ is a subset of $P$, then we will write
$\dnw X=\{y\in P\colon (\exists x\in X)(y\leq x)\}$,
$\upw X=\{y\in P\colon (\exists x\in X)(y\ge x)\}$ and say
that $X$ is a \emph{lower set} (resp. \emph{upper set})
when $X=\dnw X$ (resp. $X=\upw X$). When $X=\{a\}$, we will
sometimes write
$\dnw a$ instead of $\dnw\{a\}$. We will denote by
$\wedge$, $\bigwedge$ (resp. $\vee$, $\bigvee$) the
greatest lower bound (resp. the least upper bound) partial
operations in $P$.\smallskip

If $K$ is a convex subset of a
topological vector space, we will denote by
$\partial_\mathrm{e}K$ its extreme boundary (set of extreme
points of $K$) and by $\mathrm{Aff}(K)$ the space of all
affine continuous real-valued functions on $K$.\smallskip

Let $G$ be a \poag. An \emph{order-unit} of $G$ is an
element $u$ of $G^+$ such that
$(\forall x\in G)(\exists n\in\NN)(x\leq nu)$.
Say that $G$ is \emph{unperforated} when
it satisfies, for all $m\in\NN$, the statement
\((\forall x)(mx\geq0\Rightarrow x\geq0)\). Say that $G$ is
\emph{Archimedean} when for all elements $a$,
$b\in G$, $(\forall n\in\ZZ^+)(na\leq b)$ implies $a\leq 0$.

An \emph{interpolation group} is a \poag\ satisfying the
interpolation property.
A \emph{dimension group} is a directed, unperforated
interpolation group.

If $(G,u)$
is a \poag\ with order-unit, we will denote by $\ess(G,u)$
the state space of
$(G,u)$ (\emph{i.e.}, the set of all normalized positive
homomorphisms from $G$ to $\RR$), by
$\bphi_{(G,u)}$ the natural map from
$G$ to $\mathrm{Aff}(\ess(G,u))$ and by $\bpsi_{(G,u)}$ the
natural map from $G$ to
$\mathbf{C}(\partial_\mathrm{e}\ess(G,u),{\RR})$.

\section{Preliminaries; intervals, multiplier
groups}\label{Prel}

\begin{parag}\label{Pa:BasicDef} We shall mainly use the
notations of
\cite{Wehr96}. Thus if $(A,+,0,\leq)$ is a commutative
preordered monoid (\emph{i.e.}, $(A,+,0)$ is a commutative
monoid and
$\leq$ is a partial preordering on $A$ compatible with
$+$), we shall denote its \emph{positive cone} by
$A^+=\{x\in A\colon 0\leq x\}$ and define a new
preordering $\leq^+$ on $A$ by putting
\[ x\leq^+y\Longleftrightarrow(\exists z\geq0)(x+z=y).
\]

We shall say that $A$ is \emph{positively preordered} when
$A=A^+$. If $A$ is positively preordered, let
$\grot(A)$ be the universal group (or \emph{Grothendieck
group}) of $A$, and for all $a\in A$, denote by $[a]$ the
image of $a$ in $\grot(A)$ (thus $[a]=[b]$ if and only if
there exists
$c$ such that $a+c=b+c$); it is easy to verify that
$\grot(A)^+=\{[a]\colon a\in A\}$ is the positive cone of a
structure of partially preordered abelian group on
$\grot(A)$, which, if $A$ is positively \emph{ordered}, is
a \poag, see
\cite[Lemma 1.2]{Wehr96}.

If $A$ is a positively preordered \cm, then for all
$d\in A$, the \emph{ideal generated by $d$} is the submonoid
\[
A\res d=\{x\in A\colon (\exists n\in\NN)(x\leq^+nd)\}.
\]

Moreover, one can define a monoid congruence $\approx_d$ on
$A$ by putting
\[ x\approx_dy\Longleftrightarrow (\exists
n\in\NN)(x+nd=y+nd).
\] Note that if $x$ and $y$ are two elements of $A\res d$,
then
$x$ and $y$ have the same image in $\grot(A\res d)$ if and
only if $x\approx_dy$; thus
$\grot(A\res d)^+=(A\res d)/\!\!\approx_d$. For all
$x\in A$, denote by $[x]_d$ the equivalence class of $x$
under
$\approx_d$.
\end{parag}

\begin{parag}\label{Pa:BasicAx} We will need in this paper
five axioms among those introduced in \cite[Section
1]{Wehr96}. All the symbols used in these axioms will lie
among $+$, $\leq$ and $\approx_d$:

\begin{itemize}
\item IA (interval axiom) is
$(\forall\a,\b,\c,\d)\mathrm{IA}(\a,\b,\c,\d)$ where
$\mathrm{IA}(\a,\b,\c,\d)$ is
\[
\d\le\a+\c,\b+\c\Rightarrow(\exists\x) (\x\leq\a,\b\
\mathrm{and}\ \d\leq\x+\c).
\]

\item WIA (weak interval axiom) is
$(\forall\a,\b,\c)\mathrm{WIA}(\a,\b,\c)$ where\break
\noindent$\mathrm{WIA}(\a,\b,\c)$ is
\[
\a+\c=\b+\c\Rightarrow(\exists\x) (\x\leq\a,\b\
\mathrm{and}\
\a+\c=\x+\c).
\]

\item RD (Riesz decomposition property) is
$(\forall\a,\b,\c)\mathrm{RD}(\a,\b,\c)$ where\break
\noindent$\mathrm{RD}(\a,\b,\c)$ is
\[
\c\le\a+\b\Rightarrow(\exists\x,\y)
(\x\leq\a\ \mathrm{and}\ \y\leq\b\ \mathrm{and}\ \c=\x+\y).
\]

\item REF (refinement property) is
$(\forall\a_0,\a_1,\b_0,\b_1)\mathrm{REF}(\a_0,\a_1,\b_0,\b_1)$
where $\mathrm{REF}(\a_0,\a_1,\b_0,\b_1)$ is
\begin{multline*}
\a_0+\a_1=\b_0+\b_1\Longrightarrow
(\exists\c_{00},\c_{01},\c_{10},\c_{11})\\
(\a_0=\c_{00}+\c_{01}\ \mathrm{and}\ \a_1=\c_{10}+\c_{11}\
\mathrm{and}\ 
\b_0=\c_{00}+\c_{10}\ \mathrm{and}\ \b_1=\c_{01}+\c_{11}).
\end{multline*}

\item REF$'$ is $(\forall\d)\mathrm{REF}'(\d)$ where for
every
\cm\ $A$ and every $d\in A$, $A$ satisfies
$\mathrm{REF}'(d)$ when $\grot(A\res d)^+$ satisfies REF.
\end{itemize}
\end{parag}

Recall that the set $\La(A)$ of all
\emph{intervals} of
$(A,\leq)$ can be equipped with a natural structure of \com,
where the addition is given by 
\[
\ag+\bg=\dnw\{x+y\colon x\in\ag\
\mathrm{and}\ y\in\bg\},
\]
and the order on $\La(A)$ is just the inclusion. Note
that the positive cone of $\La(A)$ is just
$\{\ag\in\La(A)\colon 0\in\ag\}$; we will call
\emph{positive intervals} the elements of this positive
cone.

We restate here
\cite[Proposition 1.5]{Wehr96}:

\begin{proposition}\label{P:LaA,Aplus}
Let $A$ be a \com.
Then one can define two maps
\begin{gather*}
\varphi\colon \La(A)^+\to\La(A^+),\ 
\ag\mapsto\ag\cap A^+\\
\mathrm{and}\\
\psi\colon \La(A^+)\to\La(A)^+,\ \ag\mapsto\dnw\ag
\end{gather*}
and they are mutually inverse isomorphisms of
ordered monoids.\qed
\end{proposition}

In regard of this result, we will often identify positive
intervals of $A$ and intervals of $A^+$.

The following lemma is an abstract setting of
\cite[Theorem 2.7]{Good96} and the proof is essentially the
same; we write it here for convenience of the reader.

\begin{lemma}\label{L:SuffRef}
Let $A$ be a positively
\emph{ordered} monoid, let $B$ be a submonoid of $A$, and
let
$d\in B$. Suppose that the following conditions are
satisfied:

\begin{itemize}
\item[\rm (i)] $A$ satisfies \textrm{WIA};

\item[\rm (ii)] $B$ satisfies both \textrm{RD} and
\textrm{REF};

\item[\rm (iii)] For all $x\in A\res d$, there exists $y\leq
x$ in $B$ such that $x\approx_d y$.
\end{itemize}

Then $\grot(A\res d)^+$ satisfies \textrm{REF}. Thus, if in
addition $A$ is positively ordered, then
$\grot(A\res d)$ is an interpolation group.
\end{lemma}

\begin{proof} Let $x_0$, $x_1$, $y_0$, $y_1$ in
$A\res d$ such that
$[x_0]_d+[x_1]_d=[y_0]_d+[y_1]_d$. We prove that
$\grot(A\res d)^+$ satisfies
$\mathrm{REF}([x_0]_d,[x_1]_d,[y_0]_d,[y_1]_d)$. By
assumption (iii), we may assume without loss of generality
that $x_0$,
$x_1$, $y_0$, $y_1$ belong to $B$, and by definition, there
exists $n\in{\NN}$ such that
$x_0+x_1+nd=y_0+y_1+nd$. Since $A$ satisfies WIA, there
exists
$z\in A$ such that $z\leq x_0+x_1,y_0+y_1$ and
$z+nd=x_0+x_1+nd$; by assumption (iii), one may assume
without loss of generality that $z\in B$. Since $B$
satisfies RD, there exist $x'_0\leq x_0$, $x'_1\leq x_1$,
$y'_0\leq y_0$, and
$y'_1\leq y_1$ in $B$ such that $z=x'_0+x'_1=y'_0+y'_1$.
Since
$B$ satisfies REF, there exist $z_{ij}$ ($i,j<2$) in $B$
witnessing the fact that $B$ satisfies
$\mathrm{REF}(x'_0,x'_1,y'_0,y'_1)$. But 
$x_0+x_1+nd=z+nd=x'_0+x'_1+nd\leq x'_0+x_1+nd\leq
x_0+x_1+nd$ and
$x_1\in A\res d$, thus $x_0\approx_dx'_0$. Similarly, one
shows that $x_1\approx_dx'_1$ and $y_i\approx_dy'_i$ for
all
$i<2$. It follows immediately that $([z_{ij}]_d)_{i,j<2}$
witnesses
$\mathrm{REF}([x_0]_d,[x_1]_d,[y_0]_d,[y_1]_d)$ in
$\grot(A\res d)^+$. The last part of the statement results
from the fact that if $A$ is ordered, then
$(A\res d,\leq^+)$ is also ordered, thus (by
\ref{Pa:BasicDef}, or \cite[Lemma 1.2]{Wehr96})
$\grot(A\res d)$ is ordered.
\end{proof}

\begin{parag}\label{Pa:RecDefM0,M} In particular, when $G$
is a directed interpolation group, $A=\La(G^+)$ and
$\dg\in\La(G^+)$ has a countable cofinal subset, the
hypotheses above are satisfied with $B=$ the submonoid of
$A$ consisting of intervals with countable cofinal
subsets, see
\cite[Proposition 2.5 and Lemma 2.6]{Good96}; thus we
recover the statement of
\cite[Theorem 2.7]{Good96}. Let us recall the
correspondence between the definitions here and there:
\begin{gather*} M_0(G,\dg)=\La(G^+)\res\dg,\\
M(G,\dg)=\grot(M_0(G,\dg))\qquad\mbox{(multiplier group).}
\end{gather*}
\end{parag}

We shall meet in the coming sections the analogues of
$\La(G)$,
$M_0(G,\dg)$ and $M(G,\dg)$ for \emph{closed} intervals.
\bigskip

\section{Norm-closed intervals and affine \lsc\ functions
in the norm-complete case}\label{NClAff}

\begin{definition}\label{D:GapproxH} Let $(H,u)$ be a
\poag\ with order-unit, let $G$ be a subgroup of $H$, let
$a\in H$. We will say that $G$ \emph{approximates} $a$
when for all
$\varepsilon>0$, there are $n\in{\NN}$ and $x\in G$ such
that
$\|na-x\|_u\leq n\varepsilon$.
\end{definition}

Note that the definition above does not depend on the choice
of the order-unit $u$.

\begin{lemma}\label{L:barGsbgrp} In the context of
Definition~\textrm{\ref{D:GapproxH}}, the set
\[
\overline{G}=\{x\in H\colon G\ \textrm{approximates}\ x\}
\] is a subgroup of $H$ containing $G$.
\end{lemma}

\begin{proof} It is clear that $G\subseteq\overline{G}$. Let
$x$ and $y$ in $\overline{G}$ and let $\varepsilon>0$.
There are
$m$ and $n$ in ${\NN}$ and $x'$ and $y'$ in $G$ such that
$\|mx-x'\|_u\leq m\varepsilon/2$ and
$\|ny-y'\|_u\leq n\varepsilon/2$. Thus a simple calculation
yields the inequality
\(\|mn(x-y)-(nx'-my')\|\leq mn\varepsilon\) with
\(nx'-my'\in G\); thus we obtain \(x-y\in\overline{G}\).
\end{proof}

\begin{lemma}\label{L:CommDiag} In the context of
Definition~\ref{D:GapproxH} and Lemma~\ref{L:barGsbgrp},
suppose that $u\in G$. Then the restriction map
\[
\rho\colon \ess(\overline{G},u)\to\ess(G,u),\ s\mapsto
s\res_G
\] is an affine homeomorphism; therefore one has a
commutative diagram of homomorphisms of \poag s as follows:

\begin{picture}(100,100)(-130,-20)

\put(5,10){\vector(2,1){80}}
\put(90,50){\vector(0,-1){40}}
\put(10,0){\vector(1,0){55}}

\put(95,30){\makebox(0,0)[l]{$\theta$}}
\put(40,-5){\makebox(0,0)[t]{$\bphi_{(G,u)}$}}

\put(-10,0){\makebox(0,0){$(G,u)$}}
\put(110,0){\makebox(0,0) {$(\mathrm{Aff}(\ess(G,u)),1)$}}
\put(90,60){\makebox(0,0){$(\overline{G},u)$}}

\end{picture}

\noindent and if $H$ is Ar\-chi\-me\-de\-an, then $\theta$
is an embedding of ordered groups.
\end{lemma}

\begin{proof} Since $\rho$ is affine continuous, and also
surjective, see \cite[Corollary 4.3]{Good86}, it suffices
to prove that $\rho$ is one-to-one. Thus let $s$, $t$ in
$\ess(\overline{G},u)$ such that
$s\res_G=t\res_G$. Let $x\in\overline{G}$. For all
$\varepsilon>0$, there exist $n\in{\NN}$ and $y\in G$ such
that
$\|nx-y\|_u\leq n\varepsilon$; thus
$|ns(x)-s(y)|\leq n\varepsilon$ and $|nt(x)-t(y)|\leq
n\varepsilon$, thus, since $s(y)=t(y)$,
$|s(x)-t(x)|\leq 2\varepsilon$. Letting $\varepsilon$
evaporate yields $s(x)=t(x)$; whence $s=t$. Thus $\rho$ is
an affine homeomorphism. Then define $\theta$ by putting
$\theta(x)(s)=\rho^{-1}(s)(x)$; since $\rho$ is an affine
homeomorphism, $\theta$ satisfies the required properties.
The conclusion for $H$ Ar\-chi\-me\-de\-an results from
\cite[Theorem 4.14]{Good86}.
\end{proof}

\begin{lemma}\label{L:ApproxLSC} Let $(G,u)$ be a dimension
group with order-unit. Put $S=\ess(G,u)$,
$\phi=\bphi_{(G,u)}$ and
$A=\{\phi(x)/2^n\colon x\in G\ \mathrm{and}\ n\in{\NN}\}$.
Then the following properties hold:

\begin{itemize}
\item[\rm (a)] For all $q\colon S\to{\RR}$ convex \lsc, we
have $q=\bigvee\{f\in A\colon f\ll q\}$ (the supremum being
meant pointwise).

\item[\rm (b)] For all $q\colon S\to{\RR}\cup\{+\infty\}$
concave
\lsc, the set $\ddnw_\phi q=\{x\in G\colon \phi(x)\ll q\}$
is an interval of $G$.

\item[\rm (c)] For all $\ag\in\La(G)$, all
$\lambda\in{\RR}$ and all $p\colon
S\to{\RR}\cup\{-\infty\}$
\usc, if $p\ll\bigvee\phi[\ag]+\lambda$, then there exists
$a\in\ag$ such that
$p\ll\phi(a)+\lambda$.
\end{itemize}
\end{lemma}

\begin{proof} (a) Let $s\in S$ and let $\alpha<q(s)$. By
\cite[Proposition 11.8]{Good86}, there exists
$g\in\mathrm{Aff}(S)$ such that $g\ll q$ and $\alpha<g(s)$.
There exists $\varepsilon>0$ such that
$\alpha+\varepsilon<g(s)$ and $g+\varepsilon\ll q$. By
\cite[Theorem 7.9]{Good86}, there exists $f\in A$ such that
$|f-g|<\varepsilon$. Thus $f\in A$, $f\ll q$ and
$f(s)>\alpha$. Thus $q(s)=\bigvee\{f(s)\colon f\in A\
\mathrm{and}\ f\ll q\}$.
\smallskip

(b) Since $q$ is \lsc\ and $S$ is compact, $q$ is bounded
below and thus $\ddnw_\phi q\ne\emptyset$. It is trivial
that
$\ddnw_\phi q$ is a lower set. Finally, let $a$ and $b$ in
$\ddnw_\phi q$. Thus 
$p=\phi(a)\vee\phi(b)$ is a convex continuous (thus \usc)
function from $S$ to ${\RR}$ and $p\ll q$, thus, by
\cite[Theorem 11.12]{Good86}, there exists
$f\in\mathrm{Aff}(S)$ such that $p\ll f\ll q$. Let
$\varepsilon>0$ such that
$p\ll f-\varepsilon\ll f+\varepsilon\ll q$. There exists,
by
\cite[Theorem 7.9]{Good86}, $g\in A$ such that 
$|f-g|\le\varepsilon$; thus $p\ll g\ll q$. Write
$g=\phi(x)/2^n$ where $x\in G$ and $n\in{\NN}$. Since $G$ is
unperforated, we also have, by \cite[Corollary
4.13]{Good86},
$2^na,2^nb\leq x$, thus, by \cite[Proposition
2.21]{Good86}, there exists $c\in G$ such that $a,b\leq c$
and $2^nc\leq x$; it follows immediately that $\phi(c)\ll
q$.\smallskip

(c) By definition, we have $S=\bigcup\{U_a\colon a\in\ag\}$
where we put
$U_a=\{s\in S\colon p(s)<s(a)+\lambda\}$. Since the
$U_a$'s are open and $S$ is compact, there are $n\in{\NN}$
and $a_i$ ($i<n$) in $\ag$ such that
$S=\bigcup_{i<n}U_{a_i}$. Since
$\ag$ is upward directed, there exists $a\in\ag$ such that
$(\forall i<n)(a_i\leq a)$. It follows immediately that
$p\ll\phi(a)+\lambda$.
\end{proof}

Now let us recall some terminology from \cite{Alf71}. Let
$S$ be a compact convex set in a locally convex topological
vector space $E$. For every function $q\colon
S\to\RR\cup\{+\infty\}$ which is bounded below (which
happens in particular when $f$ is \lsc), one defines the
\emph{lower envelope} $\check{q}$ of
$q$ by the formula
\(\check{q}=\bigvee\mathfrak{a}(q)\) where we put
\[
\mathfrak{a}(q)=\{f\in\mathrm{Aff}(S)\colon f\leq q\}.
\] Furthermore, by \cite[Comments page 4]{Alf71},
$\check{q}$ is convex \lsc\ and we also have
\(\check{q}=\bigvee\mathfrak{b}(q)\) where we put
\[
\mathfrak{b}(q)=\{f\in\mathrm{Aff}(S)\colon f\ll q\}.
\]

\begin{lemma}\label{L:LSCgivesDir} Let $S$ be a compact
convex set in a locally convex topological vector space, let
$q\colon S\to{\RR}\cup\{+\infty\}$ be a concave \lsc\
function. Then the following assertions hold:
\begin{itemize}
\item[\rm (a)]
\(q\restriction_{\partial_{\mathrm{e}}S}=
\check{q}\restriction_{\partial_{\mathrm{e}}S}\).

\item[\rm (b)] If in addition $S$ is a Choquet simplex,
then both $\mathfrak{a}(q)$ and $\mathfrak{b}(q)$ are
upward directed and $\check{q}$ is affine.
\end{itemize}
\end{lemma}

\begin{proof} Part (a) follows from Herv\'e's Theorem
\cite[Proposition I.4.1]{Alf71}. The fact that both
$\mathfrak{a}(q)$ and $\mathfrak{b}(q)$ are upward directed
results from Edwards' Theorem (for example, to prove that
$\mathfrak{a}(q)$ is upward directed, one applies Edwards'
Theorem to
$f\vee g$ and $q$, for $f$, $g\in\mathfrak{a}(q)$). The
rest of part (b) follows from
\cite[Theorem II.3.8]{Alf71} (both applied to $-q$).
\end{proof}

For every \poag\ with order-unit $(G,u)$, denote by
$\Sigma(G,u)$ the set of all functions from
$\ess(G,u)$ to
${\RR}\cup\{+\infty\}$ of the form
$\bigvee\bphi_{(G,u)}[\ag]$ where $\ag\in\La(G)$. Thus all
elements of $\Sigma(G,u)$ are affine
\lsc\ functions from $\ess(G,u)$ to
${\RR}\cup\{+\infty\}$. As in \cite{Good96}, for every
compact convex set $K$ in a topological linear space, we
shall denote by $\varLambda(K)$ the additive monoid of all
affine
\lsc\ functions from $K$ to
${\RR}\cup\{+\infty\}$, ordered componentwise. It may of
course happen that $\Sigma(G,u)$ is a proper subset of
$\varLambda(\ess(G,u))$. We shall omit the proof of the
following lemma, which is straightforward.\smallskip

\begin{lemma}\label{L:LaIsSig} The set $\Sigma(G,u)$ is a
submonoid of $\varLambda(\ess(G,u))$, and the map
$\bigvee_\phi\colon \ag\mapsto\bigvee\phi[\ag]$ is a
homomorphism of ordered monoids from $\La(G)$ to
$\Sigma(G,u)$.\qed
\end{lemma}

\begin{lemma}\label{L:DnwPhiCl} Let $(G,u)$ be an \anog; put
$\phi=\bphi_{(G,u)}$. Let $q\in\Sigma(G,u)$. Then the set
$\dnw_\phi q=\{x\in G\colon \phi(x)\leq q\}$ is a
norm-closed interval of $G$.
\end{lemma}

\begin{proof} Put $S=\ess(G,u)$. Let $\ag\in\La(G)$ such
that
$q=\bigvee\phi[\ag]$. For all $s\in S$, the map
$G\to{\RR},\ x\mapsto\phi(x)(s)$ is continuous, thus
$\dnw_\phi q$ is norm-closed. Since $\ag\ne\emptyset$, 
$\dnw_\phi q$ is non\-emp\-ty. It is trivial that
$\dnw_\phi q$ is a lower set. Finally, let $a$,
$b\in\dnw_\phi q$. We prove first a\smallskip

\noindent{\bf Claim.} {\sl For all $n\in\ZZ^+$, there exist
$c\in\ag$ and $v\in G^+$ such that $2^nv\leq u$ and
$a,b\leq c+v$.}\smallskip

\noindent{\it Proof of Claim.} Since
$\phi(a)\vee\phi(b)\ll\bigvee\phi[\ag]+2^{-n}$, there
exists by Lemma~\ref{L:ApproxLSC} (c) an element $c$ of
$\ag$ such that
$\phi(a)\vee\phi(b)\ll\phi(c)+2^{-n}$. Therefore,
$0,2^n(a-c),2^n(b-c)\leq u$ thus, by \cite[Proposition
2.21]{Good86}, there exists $v\in G^+$ such that 
$a,b\leq c+v$ and
$2^nv\leq u$.\qed~Claim.\smallskip

In particular, there exists $c_0\in\ag$ such that
$a,b\leq c_0+u$; put $u_0=u$. Let $n\in\ZZ^+$ and suppose
having constructed $c_n\in\ag$ and $u_n\in G^+$  such that
$a,b\leq c_n+u_n$ and $2^nu_n\leq u$. By the Claim, there
exist $c\in\ag$ and $u_{n+1}\in G^+$ such that
$a,b\leq c+u_{n+1}$ and $2^{n+1}u_{n+1}\leq u$; without
loss of generality, $c_n\leq c$. Then it is easy to verify
that
$a-u_{n+1},b-u_{n+1},c_n\leq c,c_n+u_n$ thus, by
interpolation, there exists $c_{n+1}\in G$ such that
$a-u_{n+1},b-u_{n+1},c_n\leq c_{n+1}\leq c,c_n+u_n$. Since
$c_{n+1}\leq c\in\ag$, we have $c_{n+1}\in\ag$.
Furthermore,
$0\leq c_{n+1}-c_n\leq u_n$ thus $\|c_{n+1}-c_n\|_u\leq
2^{-n}$, and
$a,b\leq c_{n+1}+u_{n+1}$. Therefore the sequence
$\langle c_n\colon n\in\ZZ^+\rangle$ thus constructed is an
increasing Cauchy sequence; thus it converges to some
$c\in G$. Since
$\dnw_\phi q$ is norm-closed, we have $c\in\dnw_\phi q$.
Since $G$ is Archimedean and by
\cite[Proposition 7.17]{Good86}, 
$a,b\leq c$. Thus \(\dnw_\phi q\) is upward directed.
\end{proof}

It is very strange that the hypotheses of
Lemma~\ref{L:DnwPhiCl} cannot be weakened to arbitrary
affine
\lsc\ functions
$q\colon S\to{\RR}\cup\{+\infty\}$, even for $q$ continuous
real-valued and $G$ norm-discrete, as the following example
shows.

\begin{examplepf}\label{E:DnwNotDir} An Ar\-chi\-me\-de\-an
norm-discrete dimension group with order-unit \((G,e)\) and
\(q\in\mathrm{Aff}(\ess(G,e))^+\) such that, putting
\(\phi=\bphi_{(G,e)}\),
\(\dnw_\phi q=\{x\in G\colon \phi(x)\leq q\}\) is not upward
directed.
\end{examplepf}

\begin{proof} Put \(X=\{0,1\}^{\ZZ^+}\) be the Cantor
space, endowed with its natural product topology,
corresponding to the metric given by the formula
\noindent\(d(x,y)=\sum_{n\in{\ZZ^+}}2^{-n-1}|x(n)-y(n)|\).
Put \(E=\mathbf{LSC}^{\pm}_\mathrm{b}(X,{\ZZ}^+)\) and
\(F=\mathbf{LSC}^{\pm}_\mathrm{b}(X,{\RR}^+)\) as defined
in the Introduction. We have seen in
\cite[Proposition 3.5]{Wehr96} that $F$ is an
Ar\-chi\-me\-de\-an \poag, thus it is also the case for
$E$ (which is an ordered subgroup of $F$).\smallskip

Now, let $\alpha$ and $\beta$ be any two distinct elements
of $X$, and put $U=X\setminus\{\alpha\}$ and
$V=X\setminus\{\beta\}$; then put $a=\chi_U$, $b=\chi_V$
and $e=a+b$. Let finally $H$ (resp. $G$) be the ideal of
$F$ (resp. $E$) generated by $e$. Note that both
$a$ and $b$ (thus also $e$) belong to $G$ and that
$G\subseteq H$ (in fact $G=H\cap E$). By
\cite[Lemma 3.4 and Proposition 3.5]{Wehr96}, $E$ is an
interpolation group, thus
$G$ is also an interpolation group. By definition, $e$ is an
order-unit of $G$. For all $f\in G$ such that $\|f\|_e\leq
1/3$, we have $-e\leq 3f\leq e$, thus, since $f$ is
${\ZZ}$-valued and $e=\chi_U+\chi_V$, $f=0$; thus $G$ is
norm-discrete.\smallskip

Now, let
\(g\colon X\to[0,\,1],\
x\mapsto\mathrm{min}\{d(x,\alpha),d(x,\beta)\}\). Since
$g$ is continuous, it belongs to $F$. Since $0\leq g\leq e$
and
$g$ is continuous, we have $0\leq^+g\leq^+e$ and thus
$g\in H$. Let us prove that $G$ approximates $g$,
\emph{i.e.}, $g\in\overline{G}$ with the notation of
Lemma~\ref{L:barGsbgrp}. Thus let $\varepsilon>0$. Pick
$n\in\NN$ such that $n\varepsilon\geq1$. Since
\((ng)[X]\subseteq[0,\,n]\), we have
\(X=\bigcup_{0\leq k\leq n}(ng)^{-1}(k-1,\,k+1)\). Since
$X$ is an ultrametric space, it satisfies, by \cite[Lemma
3.4]{Wehr96}, the open reduction property \cite[Definition
3.2]{Wehr96} and thus, by
\cite[Lemma 3.3]{Wehr96}, there are clopen subsets $W_k$
($0\leq k\leq n$) of $X$ such that
\(X=\bigsqcup_{0\leq k\leq n}W_k\) and, for all
\(k\in\{0,1,\ldots,n\}\),
\(W_k\subseteq(ng)^{-1}(k-1,\,k+1)\). Put
\(h=\sum_{0\leq k\leq n}k\cdot\chi_{W_k}\). Then $h$ is
continuous and ${\ZZ}$-valued, thus $h\in G$, and 
\(-1\leq ng-h\leq 1\), thus
\emph{a fortiori} \(-e\leq ng-h\leq e\). Since $ng-h$ is
continuous, we have in fact \(-e\leq^+ng-h\leq^+e\), whence
$\|ng-h\|_e\le 1\leq n\varepsilon$. Thus
$G$ approximates $g$. By Lemma~\ref{L:barGsbgrp}, $G$ also
approximates $f=e-g$. Since $g\leq^+a,b$, we have
$a,b\leq^+f$.\smallskip

However, suppose that there exists $c\in G$ such that 
$a,b\leq^+c\leq^+f$. Put $d=e-c$. Then $d\in G$ and 
$g\leq^+d\leq^+a,b$. Since $g$ is positive, $0\leq a,b\leq
1$ and $d$ is ${\ZZ}$-valued, there exists $W\subseteq X$
such that
$d=\chi_W$. Since $d\leq a,b$, we have $W\subseteq U\cap V$.
Since $(\forall x\in U\cap V)(g(x)>0)$, we obtain $W=U\cap
V$. It follows that $\chi_{U\cap V}\leq^+\chi_U$, whence
$\chi_{\{\beta\}}$ is \lsc, a contradiction.\smallskip

So we have proved that there exists no $c\in G$ such that
$a,b\leq^+c\leq^+f$. Now, let $\theta$ be the natural
homomorphism of ordered groups from $\overline{G}$ to
$\mathrm{Aff}(\ess(G,e))$ given by Lemma~\ref{L:CommDiag}.
Since
$H$ is Ar\-chi\-me\-de\-an (it is an ideal of $F$ and $F$ is
Ar\-chi\-me\-de\-an), $\theta$ is an order-embedding. Put
$q=\theta(f)$. Then $q\in\mathrm{Aff}(\ess(G,e))$ but there
exists no
$x\in G$ such that $a,b\leq^+x$ and $\phi(x)\leq q$.
\end{proof}

Note that no \poag\ satisfying the properties of
Example~\ref{E:DnwNotDir} can be lattice-ordered. More
generally, one can easily prove, using
Lemma~\ref{L:ApproxLSC} (b), that if $(G,u)$ is a dimension
group with order-unit satisfying the
$(2,\aleph_0)$-interpolation property, then for all
$q\colon \ess(G,u)\to{\RR}\cup\{+\infty\}$ concave \lsc,
the set
$\{x\in G\colon \bphi_{(G,u)}(x)\leq q\}$ is an interval of
$G$; this holds of course in particular when $G$ is
lattice-ordered. Another case where this holds is the case
(neither more nor less general) where $(G,u)$ is an
Ar\-chi\-me\-de\-an norm-complete dimension vector space
with order-unit (this results immediately from
\cite[Corollary 15.8]{Good86}.\smallskip

Now, let us return back to the context of
Lemma~\ref{L:DnwPhiCl}: in Lemmas \ref{L:HomSigToLa} and
\ref{L:CharFinSig}, let $(G,u)$ be an
\anog. Put $S=\ess(G,u)$ and
$\phi=\bphi_{(G,u)}$.

\begin{lemma}\label{L:HomSigToLa} The map
$\dnw_\phi\colon\Sigma(G,u)\to\La(G),\ 
q\mapsto\dnw_\phi q$ is a homomorphism of ordered monoids.
\end{lemma}

\begin{proof} It is obvious that $\dnw_\phi$ is
order-preserving. Now, let $p$ and $q$ be elements of
$\Sigma(G,u)$. It is obvious that 
$\dnw_\phi p+\dnw_\phi q\subseteq\dnw_\phi(p+q)$.
Conversely, let $\ag$ and $\bg$ in $\La(G)$ such that
$p=\bigvee\phi[\ag]$ and $q=\bigvee\phi[\bg]$. Let 
$c\in\dnw_\phi(p+q)$, we shall prove that 
$c\in\dnw_\phi p+\dnw_\phi q$. We first prove the following
\smallskip

\noindent{\bf Claim.} {\sl For all $n\in\ZZ^+$, there are
$x\in\ag$, $y\in\bg$ and $v\in G^+$ such that $c\leq
x+y+v$ and
$2^nv\leq u$.}\smallskip

\noindent{\it Proof of Claim.} We have
$\phi(c)\ll\bigvee\phi[\ag+\bg]+2^{-n}$, thus, by
Lemma~\ref{L:ApproxLSC} (c), there are $a\in\ag$ and
$b\in\bg$ such that 
$\phi(c)\ll\phi(a+b)+2^{-n}$. Thus $0,2^n(c-a-b)\leq u$,
thus, by \cite[Proposition 2.21]{Good86} there exists $v\in
G^+$ such that $c\leq a+b+v$ and $2^nv\leq
u$.\qed~Claim.\smallskip

In particular for $n=0$, we obtain $a\in\ag$ and $b\in\bg$
such that $c\leq a+b+u$. Put $a_0=c-b-u$, $b_0=b$ and
$u_0=u$; we have
$a_0\in\ag$, $b_0\in\bg$, $0\leq u_0\leq u$ and
$c=a_0+b_0+u_0$. Let $n\in\ZZ^+$ and suppose that
$a_n\in\ag$,
$b_n\in\bg$ and $u_n\in G^+$ have been constructed such
that
$2^nu_n\leq u$ and $c=a_n+b_n+u_n$. By the Claim, there are
$a\in\ag$, $b\in\bg$ and $v\in G^+$ such that
$2^{n+1}v\leq u$ and $c\leq a+b+v$; in addition, we may
assume without loss of generality that $a_n\leq a$ and
$b_n\leq b$. It follows immediately that
$a_n,c-v-b\leq a,c-b_n$; thus, by interpolation, there
exists
$a_{n+1}\in G$ such that $a_n,c-v-b\leq a_{n+1}\leq
a,c-b_n$. Since $a_{n+1}\leq a$, we have $a_{n+1}\in\ag$,
and furthermore, $b_n,c-v-a_{n+1}\leq c-a_{n+1},b$, thus,
by interpolation, there exists $b_{n+1}\in G$ such that
$b_n,c-v-a_{n+1}\leq b_{n+1}\leq c-a_{n+1},b$. Since
$b_{n+1}\leq b$, we have $b_{n+1}\in\bg$. Moreover,
$a_{n+1}+b_{n+1}\leq c\leq a_{n+1}+b_{n+1}+v$, thus 
$u_{n+1}=c-(a_{n+1}+b_{n+1})$ lies between $0$ and $v$;
therefore,
$2^{n+1}u_{n+1}\leq u$. Since
$c=a_n+b_n+u_n=a_{n+1}+b_{n+1}+u_{n+1}$ and $a_n\leq
a_{n+1}$ and $b_n\leq b_{n+1}$, both $a_{n+1}-a_n$ and 
$b_{n+1}-b_n$ lie between $0$ and $u_n$.\smallskip

Therefore, the sequence $\langle a_n\colon
n\in\ZZ^+\rangle$ (resp.
$\langle b_n\colon n\in\ZZ^+\rangle$) is an increasing
Cauchy sequence of elements of $\ag$ (resp. $\bg$). If
$a=\lim_{n\to+\infty}a_n$ and $b=\lim_{n\to+\infty}b_n$,
then we obtain $a\in\dnw_\phi p$ and $b\in\dnw_\phi q$ and,
since $G$ is Archimedean and by
\cite[Proposition 7.17]{Good86},
$c=a+b$, thus proving that $c\in\dnw_\phi p+\dnw_\phi q$.
\end{proof}

We now come to the main lemma of this section; its
finiteness assumption will be removed in
Theorem~\ref{T:CharSig}.

\begin{lemma}\label{L:CharFinSig} Let \(q\colon
S\to{\RR}\) be an affine
\lsc\ function such that
\[ (\forall s\in\partial_\mathrm{e}S\ \mathrm{discrete})
(q(s)\in s[G]).
\] Then $q$ belongs to $\Sigma(G,u)$.
\end{lemma}

\begin{proof} For all $n\in\ZZ^+$, let
$\ag_n=\{x\in G\colon \phi(x)\ll q+2^{-n-1}\}$. By
Lemma~\ref{L:ApproxLSC} (b),
$\ag_n$ belongs to $\La(G)$. Thus $q_n=\bigvee\phi[\ag_n]$
belongs to $\Sigma(G,u)$. Furthermore,
$\ag_{n+1}\subseteq\ag_n$ thus
$q_{n+1}\leq q_n$, and $q_n\leq q+2^{-n-1}$ by definition of
$q_n$.\smallskip

\noindent{\bf Claim 1.} {\sl For all $n\in\ZZ^+$, one has 
$q\leq q_n$.}\smallskip

\noindent{\it Proof of Claim.} Let $s\in S$ and let
$\alpha<q(s)$. Since $q$ is affine (thus convex) \lsc, there
exists, by \cite[Proposition 11.8]{Good86},
$f\in\mathrm{Aff}(S)$ such that $f\ll q$ and $\alpha<f(s)$.
Thus, for all discrete
$t\in\partial_\mathrm{e}S$, $q(t)$ belongs both to $t[G]$
(by assumption) and to the interval
$(f(t),\,q(t)+2^{-n-1})$. Since in addition $f\ll
q+2^{-n-1}$, there exists, by
\cite[Theorem 13.5]{Good86}, $x\in G$ such that
$f\ll\phi(x)\ll q+2^{-n-1}$. Thus by definition,
$x\in\ag_n$, so that $q_n(s)\ge f(s)>\alpha$. This holds
for all
$\alpha<q(s)$, whence $q_n(s)\ge q(s)$.\qed~Claim
1.\smallskip

\noindent{\bf Claim 2.} {\sl Let $n\in\ZZ^+$ and let
$a\in\dnw_\phi q_n$. Then there exists $b\in\ag_{n+1}$ such
that
$b\leq a$ and $\|a-b\|_u\leq 2^{-n}$.}\smallskip

\noindent{\it Proof of Claim.} We have
$\phi(a)\leq q_n\leq q+2^{-n-1}\leq q_{n+1}+2^{-n-1}
\ll q_{n+1}+2^{-n}$, thus, by Lemma~\ref{L:ApproxLSC} (c),
there exists
$x\in\ag_{n+1}$ such that $\phi(a)\ll\phi(x)+2^{-n}$. Thus
$2^na\leq 2^nx+u$, thus, by \cite[Proposition
2.21]{Good86}, there exists $v\in G^+$ such that $a\leq
x+v$ and $2^nv\leq u$. Put $b=a-v$. Then $b\leq x$ thus
$b\in\ag_{n+1}$, and
$b\leq a$. Furthermore, $\|a-b\|_u=\|v\|_u\leq 2^{-n}$.
\qed~Claim 2.\smallskip

\noindent{\bf Claim 3.} {\sl The set $\ag=\dnw_\phi q$ is an
interval of $G$.}\smallskip

\noindent{\it Proof of Claim.} Since $q$ is \lsc, it is
bounded below and thus $\ag\ne\emptyset$. It is trivial
that $\ag$ is a lower set. Let $a$, $b\in\ag$. By
definition,
$a$, $b\in\ag_0$, thus, since $\ag_0$ is an interval, there
exists
$c_0\in\ag_0$ such that $a,b\leq c_0$.  Let $n\in\ZZ^+$ and
suppose having constructed $c_n\in\ag_n$ such that $a,b\leq
c_n$. By Claim 2, there exists $x\in\ag_{n+1}$ such that
$x\leq c_n$ and $\|c_n-x\|_u\leq 2^{-n}$. Since
$a$, $b\in\ag_{n+1}$ and that $\ag_{n+1}$ is an interval,
there exists $y\in\ag_{n+1}$ such that
$a,b\leq y$; since $a,b\leq c_n$, one may assume without
loss of generality (using interpolation) that $y\leq c_n$.
Since both $x$ and $y$ belong to $\ag_{n+1}$ and that
$\ag_{n+1}$ is an interval, there exists $z\in\ag_{n+1}$
such that $x,y\leq z$. Again using interpolation, there
exists
$c_{n+1}\in G$ such that $x,y\leq c_{n+1}\leq z,c_n$. Thus
$a,b\leq c_{n+1}$ and $c_{n+1}\in\ag_{n+1}$. Furthermore,
$0\leq c_n-c_{n+1}\leq c_n-x$, thus $\|c_n-c_{n+1}\|_u\leq
2^{-n}$. Therefore, the sequence $\langle c_n\colon
n\in\ZZ^+\rangle$ thus constructed is a decreasing Cauchy
sequence such that
$(\forall n\in\ZZ^+)(c_n\in\ag_n)$. Put
$c=\lim_{n\to+\infty}c_n$. Then $c\in\dnw_\phi q$ (because
for all
$n$, we have $q_n\leq q+2^{-n-1}$), and, since $G$ is
Archimedean and by
\cite[Proposition 7.17]{Good86},
$a,b\leq c$.\qed~Claim 3.\smallskip

Now, to conclude the proof, it suffices to prove that
$q=\bigvee\phi[\ag]$. It is trivial that
$q\ge\bigvee\phi[\ag]$. To prove the converse inequality,
it suffices, by Lemma~\ref{L:ApproxLSC} (a), to prove that
for all
$a\in G$ and all
$m\in{\NN}$, if we put $f=\phi(a)/2^m$, then
$f\ll q$ implies $f\le\bigvee\phi[\ag]$. Since
$\phi(a)\ll 2^mq_0=\bigvee\phi[2^m\ag_0]$, there exists by
Lemma~\ref{L:ApproxLSC} (c) an element $a_0$ of $\ag_0$ such
that
$\phi(a)\ll 2^m\phi(a_0)$, thus $a\leq 2^ma_0$. Let
$n\in\ZZ^+$ and suppose having constructed $a_n\in\ag_n$
such that $a\leq 2^ma_n$. By Claim 2, there exists
$x\in\ag_{n+1}$ such that
$x\leq a_n$ and $\|a_n-x\|_u\leq 2^{-n}$. Furthermore,
since 
$\phi(a)\ll 2^mq_{n+1}=\bigvee\phi[2^m\ag_{n+1}]$, there
exists by Lemma~\ref{L:ApproxLSC} (c) an element $y$ of
$\ag_{n+1}$ such that
$\phi(a)\ll 2^m\phi(y)$, thus $a\leq 2^my$; furthermore,
since
$a\leq 2^ma_n$, we may assume without loss of generality
that
$y\leq a_n$. Since both $x$ and $y$ belong to
$\ag_{n+1}$ and $\ag_{n+1}$ is upward directed, there
exists
$z\in\ag_{n+1}$ such that $x,y\leq z$. By interpolation,
there exists $a_{n+1}\in G$ such that 
$x,y\leq a_{n+1}\leq z,a_n$. Thus $a_{n+1}\in\ag_{n+1}$ and
$a\leq 2^ma_{n+1}$, and, in addition,
$0\leq a_n-a_{n+1}\leq a_n-x$, whence
$\|a_n-a_{n+1}\|_u\leq 2^{-n}$. Therefore, the sequence
$\langle a_n\colon n\in\ZZ^+\rangle$ thus constructed is a
decreasing Cauchy sequence such that for all
$n\in\ZZ^+$,
$a\leq 2^ma_n$ and $a_n\in\ag_n$. It follows immediately
that
$\bar a=\lim_{n\to+\infty}a_n$ belongs to $\dnw_\phi q$ and
that
$a\leq 2^m\bar a$. Hence,
$f=\phi(a)/2^m\le\phi(\bar a)\le\bigvee\phi[\ag]$. Thus 
$q\le\bigvee\phi[\ag]$, and this completes the proof.
\end{proof}

This yields a positive solution to \cite[Problem
13]{Good86}:

\begin{theorem}\label{T:SandwSig} Let $(G,u)$ be an \anog;
put
$S=\ess(G,u)$ and $\phi=\bphi_{(G,u)}$. Let 
$p\colon S\to{\RR}\cup\{-\infty\}$ be convex \usc\ and
$q\colon S\to{\RR}\cup\{+\infty\}$ be concave \lsc\ such
that
$p\leq q$ and for all discrete $s\in\partial_\mathrm{e}S$, 
$\{p(s),q(s)\}\subseteq s[G]\cup\{-\infty,+\infty\}$. Then
there exists $x\in G$ such that $p\le\phi(x)\leq q$.
\end{theorem}

Note that the answer would be the same if instead of
considering only one function $p$ and one function $q$, one
would have finitely many convex \usc\ $p_i$ ($i<m$) and
concave \lsc\ $q_j$ ($j<n$) such for all $i<m$ and $j<n$,
$p_i\leq q_j$ and for all discrete
$s\in\partial_\mathrm{e}S$,  
$\{p_i(s),q_j(s)\}\subseteq s[G]\cup\{-\infty,+\infty\}$:
it suffices to apply Theorem~\ref{T:SandwSig} to
$\bigvee_{i<m}p_i$ and 
$\bigwedge_{j<n}q_j$. In particular, if $q\ge 0$ in the
statement of Theorem~\ref{T:SandwSig}, then one can take
$x\ge 0$ --- just apply Theorem~\ref{T:SandwSig} to $0$ and
$p$ on one side,
$q$ on the other side.\smallskip

\begin{proof} Since $p$ is \usc\ and $S$ is compact, $p$ is
bounded above. Similarly, $q$ is bounded below. Therefore,
there exists $N\in{\NN}$ such that $p\leq N$ and $-N\leq q$.
Thus
$p'\leq q'$ where we put $p'=p\vee(-N)$ and $q'=q\wedge N$.
Note that $p'$ and $q'$ still satisfy the hypothesis of
Theorem~\ref{T:SandwSig}, and, in addition, they are bounded
(between $-N$ and $N$). Put
$p^*=\bigwedge\{f\in\mathrm{Aff}(S)\colon p'\leq f\}$ and
$q^*=\bigvee\{f\in\mathrm{Aff}(S)\colon f\leq q'\}$. By
Lemma~\ref{L:LSCgivesDir} (applied to $-p'$ and $q'$), we
have
$p'\leq p^*$ and $q^*\leq q'$, and $p^*$ is affine \usc\ and
$q^*$ is affine \lsc. By \cite[Theorem 11.13]{Good86}, there
exists $f\in\mathrm{Aff}(S)$ such that $p'\leq f\leq q'$;
thus
$p'\leq p^*\leq f\leq q^*\leq q'$. Furthermore, again by
Lemma~\ref{L:LSCgivesDir},
$p^*\res_{\partial_\mathrm{e}S}=p'\res_{\partial_\mathrm{e}S}$
and
$q^*\res_{\partial_\mathrm{e}S}=q'\res_{\partial_\mathrm{e}S}$. 
Therefore, $p^*$ and $q^*$ satisfy again the hypothesis of
Theorem~\ref{T:SandwSig}. But by Lemma~\ref{L:CharFinSig},
both
$q_0=q^*$ and $q_1=N-p^*$ belong to $\Sigma(G,u)$;
furthermore, 
$\phi(Nu)=N=q^*+(N-q^*)\leq q_0+q_1$, thus, by
Lemma~\ref{L:HomSigToLa}, 
$Nu\in\dnw_\phi q_0+\dnw_\phi q_1$, so that there exists 
$x\in\dnw_\phi q_0$ such that $Nu-x\in\dnw_\phi q_1$.
Therefore, $\phi(x)\leq q^*$, and $N-\phi(x)\leq N-p^*$,
\emph{i.e.}, $p^*\le\phi(x)$. It follows that one also has
$p\le\phi(x)\leq q$.
\end{proof}

On the other hand, the following counterexample shows that
the answer to the very similar \cite[Problem 19]{Good86} is
this time \emph{negative}, even for Dedekind complete
$\ell$-groups.

\begin{examplepf}\label{E:NonSandw}
Put $G=\mathbf{C}(\beta\ZZ^+,\ZZ)$ endowed with the
componentwise ordering, and let $u\in G$ be the constant
function with value $1$. Put $S=\ess(G,u)$ and
$\phi=\bphi_{(G,u)}$. Then $G$ is a Dedekind complete
$\ell$-group, but there exist an affine \usc\ function
$p\colon S\to{\RR}^+$ and an affine continuous function
$q\colon S\to{\RR}^+$ such that $p\leq q$ and
\((\forall s\in\partial_\mathrm{e}S)(p(s)\in s[G])\), but
such that there exists no \(x\in G\) such that
\(p\le\phi(x)\leq q\).
\end{examplepf}

\begin{proof}
Since $G$ is isomorphic to the additive group of all
bounded sequences of integers, it is a Dedekind
$\sigma$-complete $\ell$-group.
Put $H=\mathbf{C}(\beta\ZZ^+,\RR)$. It is
easy to see that with the terminology of
Definition~\ref{D:GapproxH},
$G$ approximates every element of $H$. Thus, by
Lemma~\ref{L:CommDiag}, the state spaces
\(\ess(G,u)\) and \(\ess(H,u)\) are isomorphic by
restriction, and, since $H$ is Archimedean, the natural map
$\theta\colon H\to\mathrm{Aff}(S)$ is an embedding of
ordered groups.

By \cite[Proposition 6.8]{Good86}, the elements of
$\ess(H,u)$ are exactly the integrals with
respect to regular Borel probability measures on
$\beta\ZZ^+$. Therefore, by
\cite[Proposition 5.24]{Good86}, the elements of
$\partial_{\mathrm{e}}\ess(H,u)$ are exactly the
evaluations at points of $\beta\ZZ^+$. By previous
paragraph, a similar statement holds for
$\partial_{\mathrm{e}}\ess(G,u)$.

Let $a_n$ ($n\in\omega$) and $b$ be the elements of
$H$ defined by the following formulas:
\begin{equation*}
a_n(\mathcal{U})=\lim_{\mathcal{U}}
\langle a_n(k)\colon k\in\ZZ^+\rangle
\ \text{and}\ 
b(\mathcal{U})=\lim_{\mathcal{U}}
\langle b(k)\colon k\in\ZZ^+\rangle
\end{equation*}
(for all $\mathcal{U}\in\beta\ZZ^+$) where we put
\[
a_n(k)=0\ \text{if}\ k<n,\quad 1\ \text{otherwise}
\]
and $b(k)=1-2^{-k}$.

Since $\langle a_n\colon
n\in\ZZ^+\rangle$ is decreasing,
$p=\bigwedge_{n\in\ZZ^+}\phi(a_n)$ is an affine \usc\
function from $S$ to ${\RR}^+$. Put $q=\theta(b)$; thus
$q\in\mathrm{Aff}(S)$. For all $n\in\ZZ^+$, we have
$2^na_n\leq 2^nb+u$; thus $p\leq q$.

Now let \(s\in\partial_\mathrm{e}S\). There exists a
ultrafilter $\mathcal{U}$ on $\ZZ^+$ such that $s$ is the
evaluation map at $\mathcal{U}$.
If $\mathcal{U}$ is principal,
\emph{i.e.}, there exists $m\in\ZZ^+$ such that
$\mathcal{U}=\{X\subseteq\ZZ^+\colon m\in X\}$, then
$s(a_n)=0$ for all
$n>m$, thus $p(s)=0$. If $\mathcal{U}$ is nonprincipal,
then $s(a_n)=1$ for all $n$, thus $p(s)=1$. Therefore, in
every case, we have $p(s)\in s[G]$.

However, suppose that there exists $x\in G$ such that
$p\le\phi(x)\leq q$. Since $\phi(x)\leq q$, we have
$x\leq b$; but $x$ is $\ZZ$-valued, thus $x\leq 0$; thus
$p\leq 0$. But if
$s$ is the limit operation with respect to a nonprincipal
ultrafilter, then $p(s)=1$, a contradiction.
\end{proof}

We now turn to positive applications of
Theorem~\ref{T:SandwSig}. First, it allows us to
characterize completely the elements of
$\Sigma(G,u)$ (thus strengthening
Lemma~\ref{L:CharFinSig}):

\begin{theorem}\label{T:CharSig} Let $(G,u)$ be an \anog;
put $S=\ess(G,u)$. Let $q\colon S\to{\RR}\cup\{+\infty\}$.
Then the following are equivalent:

\begin{itemize}
\item[\rm (i)] $q\in\Sigma(G,u)$;

\item[\rm (ii)] $q$ is affine \lsc\ and
$(\forall s\in\partial_\mathrm{e}S\ \mathrm{discrete})
\bigl(q(s)\in s[G]\cup\{+\infty\}\bigr)$.
\end{itemize}
\end{theorem}

\begin{proof} We prove the non-trivial direction. Thus let
$q$ satisfying condition (ii). Put
$\phi=\bphi_{(G,u)}$. Put $\ag=\dnw_\phi q$. We prove that
$\ag\in\La(G)$ and $q=\bigvee\phi[\ag]$. Since $q$ is \lsc,
$\ag$ is a non\-emp\-ty lower subset of $G$. Let
$a$, $b\in\ag$. Then
$\phi(a)\vee\phi(b)\leq q$; it is easy to verify that the
conditions of Theorem~\ref{T:SandwSig} are fulfilled, thus
there exists
$c\in G$ such that  $\phi(a)\vee\phi(b)\le\phi(c)\leq q$.
Since $G$ is Archimedean and by
\cite[Theorem 7.7]{Good86}, we have $a,b\leq c$.
This proves that $\ag\in\La(G)$. It is trivial that
$\bigvee\phi[\ag]\leq q$. To prove the converse
inequality, it suffices, by Lemma~\ref{L:ApproxLSC} (a),
to prove that for all $f\in\mathrm{Aff}(S)$ such that
$f\ll q$, we have $f\ll\bigvee\phi[\ag]$. Since $f$ is
bounded above, we have
$f\ll q\wedge N$ for some $N\in{\NN}$.
Let $q^*$ be the lower
envelope of $q\wedge N$ (the definition of the lower
envelope is recalled before Lemma~\ref{L:LSCgivesDir}).
Since $f\ll q\wedge N$, $S$ is compact, $f$ is
continuous and $q\wedge N$ is \lsc, there exists
$\varepsilon>0$ such that
$f+\varepsilon\leq q\wedge N$. Then it follows from the
definition of $q^*$ that $f+\varepsilon\leq q^*$.
Furthermore, by Lemma~\ref{L:LSCgivesDir},
$q^*$ is affine \lsc\ and for all
$s\in\partial_\mathrm{e}S$,
$q^*(s)=\mathrm{min}\{q(s),N\}\in s[G]$. Since $q^*$ is
bounded, it results from Lemma~\ref{L:CharFinSig} that
$q^*\in\Sigma(G,u)$, so that there exists
$\ag^*\in\La(G)$ such that $q^*=\bigvee\phi[\ag^*]$. Since
$f\ll q^*$, it results from Lemma~\ref{L:ApproxLSC} (c)
that there exists
$x\in\ag^*$ such that $f\ll\phi(x)$. Since $q^*\leq q$, we
also have $x\in\ag$. Thus 
$f\ll\phi(x)\le\bigvee\phi[\ag]$, which concludes the
proof.
\end{proof}

Now, equip ${\RR}\cup\{+\infty\}$ with the metric $d$
defined by $d(x,y)=\break\mathrm{min}\{|x-y|,1\}$ when
both $x$ and $y$ are real, and $d(x,+\infty)=1$ when $x$ is
real. Then the following corollary is a straightforward
consequence of Theorem~\ref{T:CharSig}:

\begin{corollary}\label{C:ClosSig} Let $(G,u)$ be an \anog.
Then $\Sigma(G,u)$ is closed under uniform limit in
$(\RR\cup\{+\infty\})^S$.\qed
\end{corollary}

The analogue of this result for the metric on
${\RR}\cup\{+\infty\}$ inherited from the natural metric on
$[-\infty,\,+\infty]$ is false (the sequence
$\langle -n\colon n\in\ZZ^+\rangle$ converges uniformly to
$-\infty$ in the space $[-\infty,\,+\infty]^S$ but does not
converge for the metric above to any element of
$\Sigma(G,u)$), but true for sequences which are uniformly
bounded below.

\begin{proposition}\label{P:ClosInt}
Let $(G,u)$ be a dimension group
with order-unit. Let $\ag\in\La(G)$,
let $a\in G$. Then the following are equivalent:

\begin{itemize}
\item[\rm (i)] $a$ belongs to the norm-closure
$\mathrm{Cl}(\ag)$ of $\ag$;

\item[\rm (ii)] $\phi(a)\le\bigvee\phi[\ag]$;

\item[\rm (iii)] There exists an increasing sequence of
elements of $\ag$ which norm-converges to $a$.
\end{itemize}
\end{proposition}

\begin{proof} (i)$\Rightarrow$(ii) is easy.\smallskip

(ii)$\Rightarrow$(iii) Assume (ii). We start with the
following\smallskip

\noindent{\bf Claim.} {\sl For all $n\in\ZZ^+$, there are
$x\in\ag$ and $v\in G^+$ such that $a\leq x+v$ and
$2^nv\leq u$.}\smallskip

\noindent{\it Proof of Claim.} Since
$\phi(a)\ll\bigvee\phi[\ag]+2^{-n}$, there exists by
Lemma~\ref{L:ApproxLSC} (c) an element $x$ of $\ag$ such
that
$\phi(a)\ll\phi(x)+2^{-n}$. Thus, by 
\cite[Corollary 4.13]{Good86}, $2^na\ll 2^nx+u$; thus, by
\cite[Proposition 2.21]{Good86}, there exists $v\in G^+$
such that $2^nv\leq u$ and $a\leq
x+v$.\qed~Claim.\smallskip

In particular for $n=0$, we obtain $a_0\in\ag$ and $u_0=u$
such that $a=a_0+u_0$. Let $n\in\ZZ^+$ and suppose having
constructed $a_n\in\ag$ and
$u_n\in G^+$ such that $a=a_n+u_n$ and $2^nu_n\leq u$. By
the Claim, there are $x\in\ag$ and $v\in G^+$ such that
$a\leq x+v$ and $2^{n+1}v\leq u$. Since $\ag$ is an
interval, one may assume without loss of generality that
$a_n\leq x$. By interpolation, there exists $a_{n+1}\in G$
such that
$a-v,a_n\leq a_{n+1}\leq a,x$. Put $u_{n+1}=a-a_{n+1}$;
since
$0\leq u_{n+1}\leq v$, we have $2^{n+1}u_{n+1}\leq u$.
Since
$a_{n+1}\leq x$ we have $a_{n+1}\in\ag$. Furthermore, 
$0\leq a_{n+1}-a_n\leq u_n$ thus $\|a_{n+1}-a_n\|_u\leq
2^{-n}$. Therefore, the sequence $\langle a_n\colon
n\in\ZZ^+\rangle$ thus constructed is an increasing Cauchy
sequence of elements of
$\ag$, with limit
$a$.\smallskip

(iii)$\Rightarrow$(i) is trivial.
\end{proof}

\section{The monoid of norm-closed intervals}\label{MonNCl}

In this section, we shall apply Theorem~\ref{T:SandwSig} to
a more complete study of norm-closed intervals of
\anogs. In \ref{C:LaClMon} -- \ref{L:LaSigClRef}, let
$(G,u)$ be an \anog, and put $S=\ess(G,u)$ and
$\phi=\bphi_{(G,u)}$. From now on, denote by
$\La_\mathrm{cl}(G)$ the space of all
\emph{norm-closed} intervals of $G$. From Lemmas
\ref{L:LaIsSig}, \ref{L:HomSigToLa} and
Proposition~\ref{P:ClosInt}, we deduce immediately the
following corollaries:

\begin{corollary}\label{C:LaClMon} The set
$\La_\mathrm{cl}(G)$ is closed under addition of intervals,
and the closure map
$\ag\mapsto\mathrm{Cl}(\ag)$ is a retraction from $\La(G)$
onto
$\La_\mathrm{cl}(G)$.\qed
\end{corollary}

Recall that $\Sigma(G,u)$ is an
ordered submonoid of $\varLambda(\ess(G,u))$ (see
Lemma~\ref{L:HomSigToLa}).

\begin{corollary}\label{C:SigIsLaCl} The map $\dnw_\phi$
determines an isomorphism from $\Sigma(G,u)$ onto
$\La_\mathrm{cl}(G)$, and its inverse is the map
$\bigvee_\phi$.\qed
\end{corollary}

By analogy with \cite{Good96}, for every norm-closed
positive interval (see the comments preceding
Proposition~\ref{P:LaA,Aplus}) $\dg$ of $G$, we shall put
$M_{0,\mathrm{cl}}(G,\dg)=\La_\mathrm{cl}(G)^+\res\dg$ and
$M_\mathrm{cl}(G,\dg)=\grot(\La_\mathrm{cl}(G)^+\res\dg)$.
It is to be noted that, by Corollary~\ref{C:LaClMon},
$M_{0,\mathrm{cl}}(G,\dg)$ (resp. $M_\mathrm{cl}(G,\dg)$)
is a \emph{retract} of $M_0(G,\dg)$ (resp.
$M(G,\dg)$).\medskip

\begin{parag} As shown in \cite[Theorem 3.8]{Wehr96}, there
are cases where $(G,u)$ is norm-discrete (thus
$\La(G)=\La_\mathrm{cl}(G)$) although $\La(G)^+$
satisfies  a strong negation of both REF and REF$'$
(denoted there by NR). Thus, in order to obtain positive
results, we shall focus attention on those ``countably
generated" elements of
$\La_\mathrm{cl}(G)$. The corresponding theory bears close
similarities with \cite[Section 2]{Good96}.\smallskip

Denote by $\La_{(\sigma)}(G)$ the submonoid of $\La(G)$
whose elements are those intervals of $G$ having a
countable cofinal subset; note that if such an interval is
positive, then it has a countable cofinal subset in $G^+$.
Say that an element of
$\La_\mathrm{cl}(G)$ (resp. $\Sigma(G,u)$) is
\emph{separable} when it is the image under Cl (resp.
$\bigvee_\phi$) of an element of
$\La_{(\sigma)}(G)$ (in the case of an interval, this is
strictly weaker than having a countable dense subset), and
denote by
$\La_{\sigma,\mathrm{cl}}(G)$ (resp.
$\Sigma_{\sigma}(G,u)$) the set of all separable elements
of $\La_\mathrm{cl}(G)$ (resp.
$\Sigma(G,u)$). An important difference with the case
without any cardinality restriction is that
\emph{$\La_{\sigma,\mathrm{cl}}(G)$ is no longer a retract
of
$\La_{(\sigma)}(G)$ (it may not even be a subset of
it)}.\smallskip

The following lemma is a version for norm-closed intervals
of
\cite[Lemma 2.6]{Good96}, and its proof uses this result.
\end{parag}

\begin{lemma}\label{L:PrepRefLaCl}
\begin{itemize}
\item[\rm (a)] Let
$\ag$, $\bg\in\La_\mathrm{cl}(G)^+$ and let
$\cg\in\La_{\sigma,\mathrm{cl}}(G)^+$ such that
$\cg\subseteq\ag+\bg$. Then there are $\ag'\subseteq\ag$
and
$\bg'\subseteq\bg$ in $\La_{\sigma,\mathrm{cl}}(G)^+$ such
that
$\cg=\ag'+\bg'$.

\item[\rm (b)] Let $\dg\in\La_{\sigma,\mathrm{cl}}(G)^+$.
For all
$\ag\in M_{0,\mathrm{cl}}(G,\dg)$, there exists
$\ag'\subseteq\ag$ in
$\La_{\sigma,\mathrm{cl}}(G)^+$ such that
$\ag\approx_{\dg}\ag'$.
\end{itemize}
\end{lemma}

\begin{proof} (a) Let $\cg_0\in\La_{(\sigma)}(G)^+$ such
that
$\cg=\mathrm{Cl}(\cg_0)$. Thus
$\cg_0\subseteq\ag+\bg$, thus, by \cite[Lemma 2.6]{Good96},
there are
$\ag_0\subseteq\ag$ and
$\bg_0\subseteq\bg$ in $\La_{(\sigma)}(G)^+$ such that
$\cg_0=\ag_0+\bg_0$. Take $\ag'=\mathrm{Cl}(\ag_0)$ and
$\bg'=\mathrm{Cl}(\bg_0)$.\smallskip

(b) There exist $n\in{\NN}$ and
$\bg\in\La_\mathrm{cl}(G)^+$ such that $\ag+\bg=n\dg$. By
(a) there are $\ag'\subseteq\ag$ and
$\bg'\subseteq\bg$ in $\La_{\sigma,\mathrm{cl}}(G)^+$ such
that $n\dg=\ag'+\bg'$. Then $\ag'$ satisfies the required
conditions.
\end{proof}

\begin{lemma}\label{L:LaSigClRef} The monoid 
$\La_{\sigma,\mathrm{cl}}(G)^+$ satisfies the refinement
property.
\end{lemma}

\begin{proof} Let $\ag_0$, $\ag_1$, $\bg_0$, $\bg_1$ in
$\La_{\sigma,\mathrm{cl}}(G)^+$ such that
$\ag_0+\ag_1=\bg_0+\bg_1$. For all $i<2$, let $\ag'_i$
(resp.
$\bg'_i$) an element of
$\La_{(\sigma)}(G)^+$ of closure $\ag_i$ (resp. $\bg_i$) and
let $\langle a_{in}\colon n\in\ZZ^+\rangle$ (resp.
$\langle b_{in}\colon n\in\ZZ^+\rangle$) be an increasing
sequence of elements of
$G^+$ which is cofinal in
$\ag'_i$ (resp. $\bg'_i$). Put $a^*_{i0}=a_{i0}$ and
$b^*_{i0}=b_{i0}$. For $n\in\ZZ^+$ even, suppose having
constructed, for all $i<2$, $a^*_{in}\in\ag_i$. There are
$b^*_{0,n+1}\in\bg_0$ and $b^*_{1,n+1}\in\bg_1$ such that
$a^*_{0n}+a^*_{1n}\leq b^*_{0,n+1}+b^*_{1,n+1}$;
furthermore, one can assume without loss of generality that
for all $i<2$, we have
$b^*_{in},b_{i,n+1}\leq b^*_{i,n+1}$. Similarly, for
$n\in\ZZ^+$ odd, if $b^*_{0n}\in\bg_0$ and
$b^*_{1n}\in\bg_1$ have been constructed, then there are
$a^*_{0,n+1}\in\ag_0$ and $a^*_{1,n+1}\in\ag_1$ such that
$b^*_{0n}+b^*_{1n}\leq a^*_{0,n+1}+a^*_{1,n+1}$ and for all
$i<2$, $a^*_{in},a_{i,n+1}\leq a^*_{i,n+1}$.

For all $i<2$, let $\ag^*_i$ (resp. $\bg^*_i$) be the
interval of $G$ generated by $\{a^*_{in}\colon n\in\ZZ^+\}$
(resp.
$\{b^*_{in}\colon n\in\ZZ^+\}$) -- thus all these intervals
belong to $\La_{(\sigma)}(G)^+$. By construction,
$\mathrm{Cl}(\ag^*_i)=\ag_i$ and
$\mathrm{Cl}(\bg^*_i)=\bg_i$, and
$\ag^*_0+\ag^*_1=\bg^*_0+\bg^*_1$. Applying the fact that 
$\La_{(\sigma)}(G)^+$ satisfies REF, see
\cite[Proposition 2.5]{Good96}, and taking closures (use
Corollary~\ref{C:LaClMon}) yields immediately
$\mathrm{REF}(\ag_0,\ag_1,\bg_0,\bg_1)$ in
$\La_{\sigma,\mathrm{cl}}(G)^+$.
\end{proof}

This allows us to deduce the following

\begin{proposition}\label{P:MClDimGrp} Let $(G,u)$ be an
\anog, let $\dg$ be the closure of a positive interval of
$G$ with a countable cofinal subset. Then
$M_\mathrm{cl}(G,\dg)$ is a dimension group.
\end{proposition}

\begin{proof} Since, by Corollary~\ref{C:LaClMon},
$\La_\mathrm{cl}(G)$ is an ordered submonoid of $\La(G)$ and
since
$M(G,\dg)$ (see \ref{Pa:RecDefM0,M}) is unperforated, see
\cite[Corollary 2.4]{Good96}, $M_\mathrm{cl}(G,\dg)$ is
also unperforated. Since
$\La_\mathrm{cl}(G)^+$ is a retract of $\La(G)^+$
(Corollary~\ref{C:LaClMon}) and that $\La(G)^+$ satisfies
IA, see \cite[Lemma 1.7]{Wehr96},
$\La_\mathrm{cl}(G)^+$ satisfies IA, thus WIA. The rest
results from Lemmas \ref{L:PrepRefLaCl},
\ref{L:LaSigClRef} and
\ref{L:SuffRef} (for
$A=\La_\mathrm{cl}(G)^+$ and
$B=\La_{\sigma,\mathrm{cl}}(G)^+$).
\end{proof}

\begin{corollary}\label{C:MClDimGrp}
Let $(G,u)$ be an
\anog\ such that $\ess(G,u)$ is metrizable. Then for all
$\dg\in\La_\mathrm{cl}(G)^+$, $M_\mathrm{cl}(G,\dg)$ is a 
dimension group.
\end{corollary}

\begin{proof}
Since $(G,u)$ is an \anog, it embeds as an ordered group in
$\mathbf{C}(\ess(G,u),\RR)$ (see
\cite[Theorem 7.7 (a)]{Good86}), and this embedding
preserves the norm. Since $\ess(G,u)$ is compact metrizable,
$\mathbf{C}(\ess(G,u),\RR)$ is separable (see for example
\cite[Proposition 5.23]{Good86}). Therefore, $(G,u)$,
endowed with its natural norm, is metrizable separable.
Thus, $\dg$ is also separable. Let $\langle a_n\colon
n\in\NN\rangle$ be a dense sequence of $\dg$. Since $\dg$
is a positive interval, there exists an increasing sequence
$\langle b_n\colon n\in\NN\rangle$ of elements of $\dg\cap
G^+$ such that $a_n\leq b_n$ for all $n$. Let $\dg'$ be the
interval generated by $\{b_n\colon n\in\NN\}$. Then $\dg$
is the closure of $\dg'$. We conclude by
Proposition~\ref{P:MClDimGrp}.
\end{proof}

In the case where $(G,u)$ is an \anog\ such that
$\partial_\mathrm{e}\ess(G,u)$ is
\emph{compact} (\emph{i.e.}, $G$ is a $\ell$-group by
\cite[Corollary 15.10]{Good86}), then the norm-closed
intervals of $G$ let themselves be described in a somewhat
more wieldy way than in Corollaries \ref{C:LaClMon} and
\ref{C:SigIsLaCl}. Indeed, let
$\psi=\bpsi_{(G,u)}$ be the natural map from $G$ to
$\mathbf{C}(\partial_\mathrm{e}\ess(G,u),{\RR})$ and let
$\Sigma_\mathrm{e}(G,u)$ be the set of all functions from
$\partial_\mathrm{e}\ess(G,u)$ to
${\RR}\cup\{+\infty\}$ of the form $\bigvee\psi[\ag]$ where
$\ag\in\La(G)$. One can then prove the following
proposition:

\begin{proposition}\label{P:SigIsSige} Let $(G,u)$ be an
Ar\-chi\-me\-de\-an norm-complete $\ell$-group with
order-unit. Put
$\phi=\bphi_{(G,u)}$ and $\psi=\bpsi_{(G,u)}$. Then one can
define an isomorphism of ordered monoids from $\Sigma(G,u)$
to $\Sigma_\mathrm{e}(G,u)$ which for all
$\ag\in\La(G)$ sends $\bigvee\phi[\ag]$ to
$\bigvee\psi[\ag]$.
\end{proposition}

\begin{proof} Clearly, it suffices to prove that for all
$\ag$ and $\bg$ in $\La(G)$, one has 
$\bigvee\phi[\ag]\le\bigvee\phi[\bg]\Leftrightarrow
\bigvee\psi[\ag]\le\bigvee\psi[\bg]$; furthermore, since
$\bigvee\phi[\cg]=\bigvee\phi[\mathrm{Cl}(\cg)]$ and
$\bigvee\psi[\cg]=\bigvee\psi[\mathrm{Cl}(\cg)]$ for all
$\cg\in\La(G)$, it suffices to prove it for $\ag$ and $\bg$ 
\emph{norm-closed}. If
$\bigvee\phi[\ag]\le\bigvee\phi[\bg]$, then, by 
Proposition~\ref{P:ClosInt}, $\ag\subseteq\bg$ thus
$\bigvee\psi[\ag]\le\bigvee\psi[\bg]$. Conversely, suppose
that  $\bigvee\psi[\ag]\le\bigvee\psi[\bg]$. Put
$X=\partial_\mathrm{e}\ess(G,u)$; thus $X$ is compact
Hausdorff. Let $a\in\ag$; we have
$\psi(a)\ll\bigvee\psi[\bg]+2^{-n}$, thus, using
compactness of
$X$, there exists $b\in\bg$ such that
$\psi(a)\ll\psi(b)+2^{-n}$. Thus $2^n(a-b)\leq u$, thus
there exists $v\in G^+$ such that $a\leq b+v$ and
$2^nv\leq u$. Put
$x=a-v$. Then $x\in\bg$ and $\|a-x\|_u\leq 2^{-n}$: this
proves that $a\in\mathrm{Cl}(\bg)$. Since $\bg$ is
norm-closed, we obtain
$\ag\subseteq\bg$, whence
$\bigvee\phi[\ag]\le\bigvee\phi[\bg]$.
\end{proof}

This allows us to construct the following example (note the
similarity with \cite[Example 7.6]{Good96}).

\begin{examplepf}\label{E:MClNotArch} Put
\(G=\mathbf{C}([0,\,1],{\RR})\), equipped with the constant
function $u$ with value $1$ as an order-unit. Then there
exists a norm-closed positive interval $\dg$ of $G$ such
that 
\(M_\mathrm{cl}(G,\dg)\) is not Ar\-chi\-me\-de\-an.
\end{examplepf}

\begin{proof} Put
\(E=\mathbf{LSC}([0,\,1],{\RR}\cup\{+\infty\})\). By 
\cite[Corollary 15.8]{Good86} and both
Corollary~\ref{C:SigIsLaCl} and
Proposition~\ref{P:SigIsSige},
\(\La_\mathrm{cl}(G)\) is isomorphic to \(E\); thus we will
argue in \(E\). For all real \(\alpha>0\),
let
\(f_{\alpha}\) be the function from \([0,\,1]\) to
\({\RR}\) defined by
\(f_{\alpha}(0)=\alpha\) and for all \(t\in(0,\,1]\), 
\(f_{\alpha}(t)=1/t\). It is easy to verify that
\(f_{\alpha}\in E^+\). Put \(d=f_{1}\); we shall prove
that
\(\grot(E^+\res d)\) is not Ar\-chi\-me\-de\-an. Put
\(a=f_{0}\). Then
\(a+f_{2}=2d\) thus \(a\in E^+\res d\), and for all
\(n\in\NN\),
\(na+f_{n+1}=(n+1)d\) thus \(n[a]\le(n+1)[d]\) in
\(\grot(E^+\res d)\). However, suppose that \([a]\le[d]\)
in
\(\grot(E^+\res d)\). Then there exist \(g\in E^+\) and
\(n\in{\NN}\) such that
\(a+g+nd=(n+1)d\), whence \(a+g=d\) since \(d\) assumes
only finite values; therefore,  \(g=d-a=\chi_{\{0\}}\) is
\lsc, a contradiction.
\end{proof}

In the example above, \(\dg\) is an unbounded interval
(although the corresponding \(d\in\Sigma_\mathrm{e}(G,u)\)
takes only finite values). We shall conclude this section by
proving that when
\(\dg\) is bounded, then \(M_\mathrm{cl}(G,\dg)\) is always
Ar\-chi\-me\-de\-an norm-complete, even though by the
results of \cite[Section 3]{Wehr96}, it may not have
interpolation.

\begin{theorem}\label{T:MClArch} Let \((G,u)\) be an
\anog, let
\(\dg\) be a bounded positive norm-closed interval of \(G\).
Then
\(M_{0,\mathrm{cl}}(G,\dg)\) is cancellative and
\(M_\mathrm{cl}(G,\dg)\) is Ar\-chi\-me\-de\-an and
norm-complete.
\end{theorem}

\begin{proof} Put as usual \(S=\ess(G,u)\) and
\(\phi=\bphi_{(G,u)}\). Put \(d=\bigvee\phi[\dg]\). Then,
by Corollary~\ref{C:SigIsLaCl},
\(M_{0,\mathrm{cl}}(G,\dg)\) is isomorphic to
\(\Sigma(G,u)^+\res d\); since $d$ is real-valued, 
\(M_{0,\mathrm{cl}}(G,\dg)\) is cancellative. Now let
$f$, $g$ and $h$ in \(\Sigma(G,u)^+\res d\) such that for
all \(n\in{\NN}\),
\(nf\leq^+ng+h\). Thus for all \(n\in{\NN}\), the map
\(h_n=g-f+(1/n)h\) is positive \lsc; since \(h\) is bounded,
\(\langle h_n\colon n\in{\NN}\rangle\) converges uniformly
to
\(g-f\), thus
\(g-f\) is positive \lsc. For all discrete
\(s\in\partial_\mathrm{e}S\),
\((g-f)(s)\in s[G]\), thus \(g-f\in\Sigma(G,u)^+\) by
Theorem~\ref{T:CharSig}; whence \(f\leq^+g\). This proves
that
\(M_\mathrm{cl}(G,\dg)\) is Ar\-chi\-me\-de\-an.\smallskip

We finally prove norm-completeness. It suffices to prove
that if \(\langle f_n\colon n\in\ZZ^+\rangle\) is a
sequence of elements of
\(\Sigma(G,d)^+\res d\) such that for all \(n\),
\(\|f_{n+1}-f_n\|_d<2^{-n-1}\), then it is convergent for
\(\|\underline{\phantom{A}}\|_d\). First, since \(d\) is
bounded, 
\(\langle f_n\colon n\in\ZZ^+\rangle\) is a Cauchy sequence
for the norm of the uniform convergence, thus it converges
uniformly to some
\(f\colon S\to\RR\); by Corollary~\ref{C:ClosSig}, $f$
belongs to
\(\Sigma(G,u)\). Furthermore, let \(N\in\NN\) such that
for all
\(n\in\ZZ^+\), 
\(\|f_n\|_d<N\). Then for all \(n\in\ZZ^+\),
\(f'_n=Nd-f_n\) belongs to \(\Sigma(G,u)^+\) and
\(\langle f'_n\colon n\in\ZZ^+\rangle\) is a Cauchy
sequence (for either norm). Thus, again by
Corollary~\ref{C:ClosSig}, it converges uniformly to some
\(f'\in\Sigma(G,u)^+\). Since \(f_n+f'_n=Nd\) for all
\(n\), we obtain that \(f+f'=Nd\); whence
\(f\in\Sigma(G,u)^+\res d\).

For all \(n\in\ZZ^+\), we have
\(\|2^{n+1}f_{n+1}-2^{n+1}f_n\|_d<1\), thus there exists
\(g_n\in\Sigma(G,u)^+\) such that
\(2^{n+1}f_n+g_n=2^{n+1}f_{n+1}+d\), \emph{i.e.},
\(f_n+g_n/2^{n+1}=f_{n+1}+d/2^{n+1}\). It follows easily
that for all
\(k\in\ZZ^+\), we have
\[
f_n+\sum_{i<k}\frac{g_{n+i}}{2^{n+i+1}}=
f_{n+k}+\sum_{i<k}\frac{d}{2^{n+i+1}},
\]
thus, letting $k$ go to infinity,
\[
f_n+\frac{g'_n}{2^n}=f+\frac{d}{2^n}\ \mathrm{where}\ 
g'_n=\sum_{i\in\ZZ^+}\frac{g_{n+i}}{2^{i+1}}.
\]
Thus \(g'_n\) is positive affine \lsc, and since \(f_n\),
\(f\) and
\(d\) lie in \(\Sigma(G,u)\), \(g'_n(s)\in s[G]\) for all
discrete
\(s\in\partial_\mathrm{e}S\); thus, by
Theorem~\ref{T:CharSig},
\(g'_n\in\Sigma(G,u)^+\). Therefore, \(2^n(f_n-f)\leq^+d\)
in
\(\Sigma(G,u)\). One can prove similarly that
\(2^n(f-f_n)\leq^+d\) in \(\Sigma(G,u)\). It follows that
\(\|f_n-f\|_d\leq2^{-n}\), so that
\(f=\lim_{n\to+\infty}f_n\) for
\(\|\underline{\phantom{A}}\|_d\). The conclusion follows.
\end{proof}

\begin{problem}\label{Pb:SpSentBarG,G} Say as in
\cite{Wehr95} that a \emph{special} sentence is a sentence
of the form
\((\forall\vec x)(\varphi\Rightarrow(\exists\vec y)\psi)\)
where \(\varphi\) and \(\psi\) are conjunctions of atomic
formulas. Is the set of all \emph{special} sentences which
are true in all structures
\((\overline{G},G,+,0,\le)\), where \((G,u)\) is an \anog\
and
\(\overline{G}=\mathrm{Aff}(\ess(G,u))\) (\(G\) being
identified with its natural image into \(\overline{G}\))
\emph{decidable}? Note that there are non-trivial
sentences to decide, as, \emph{e.g.}, the one leading to
the relatively complicated Example~\ref{E:DnwNotDir}.
\end{problem}

\begin{problem}\label{Pb:LaSigSatSD} Let \((G,u)\) be an
\anog. Do
\(\La_{(\sigma)}(G)^+\) and
\(\La_{\sigma,\mathrm{cl}}(G)^+\) always satisfy the axiom
SD of \cite{Wehr96},
\emph{i.e.}, for all \(a_0\), \(a_1\), \(b\) and \(c\) such 
that \(a_0+a_1+c=b+c\), do there exist \(b_0\), \(b_1\),
\(c_0\) and \(c_1\) such that \(b_0+b_1=b\) and
\(c_0+c_1=c\) and \(a_i+c_i=b_i+c_i\) for all \(i<2\)?
\end{problem}

\begin{ackn}
The author wishes to thank the referee for many helpful
comments.
\end{ackn}

\end{document}